\documentclass{article}

\usepackage{amsmath}
\usepackage{amssymb}
\usepackage{latexsym}
\usepackage{amsbsy}

\newcommand{\R}{\mathbb R}
\newcommand{\C}{\mathbb C}
\newcommand{\K}{\mathbb K}
\newcommand{\G}{\ensuremath{\mathcal{G}}}
\newcommand{\gs}{\ensuremath{\mathcal{G}}}
\newcommand{\esm}{\ensuremath{{\mathcal E}_m} }
\newcommand{\ns}{\ensuremath{{\mathcal N}} }
\newcommand{\beas}{\begin{eqnarray*}}
\newcommand{\eeas}{\end{eqnarray*}}
\newcommand{\cinfty}{{\mathcal C}^\infty}
\newcommand{\eps}{\varepsilon}
\newcommand{\comp}{\subset\subset}

\begin{document}

\title{Recent Progress in Special Colombeau Algebras:\\ 
Geometry, Topology, and Algebra}

\author{M. Kunzinger\footnote{Faculty of Mathematics, University of Vienna,
         Nordbergstr.\ 15, A-1090 Wien, Austria,
         Electronic mail: michael.kunzinger@univie.ac.at. Work supported by FWF-project
         P-16742 and START-project Y-237.}}
\date{}
\maketitle

\begin{abstract}Over the past few years there has been considerable progress in the
structural understanding of special Colombeau algebras. We present some of the main
trends in this development: non-smooth differential geometry, locally convex theory 
of modules over the ring of generalized numbers, and algebraic aspects of Colombeau
theory. Some open problems are given and directions of further research are outlined.

\end{abstract}

\section*{1. Introduction.} 
Colombeau algebras of generalized functions~\cite{c1, c2, MObook} are differential algebras 
that contain the vector space of Schwartz distributions as a linear subspace, 
and the space of smooth functions as 
a faithful subalgebra. Initially discovered in the context of infinite-dimensional calculus in locally 
convex spaces, such algebras have turned out to be a powerful tool in the study of singular
problems that involve 
differentiation combined with non-linear operations. In particular, Colombeau algebras quickly found 
(and continue to find) applications in the field of non-linear partial differential equations 
(e.g., \cite{MObook,Biag,NP}), 
where the application of classical distributional methods is limited by the 
impossibility to consistently define an intrinsic product of distributions~\cite{Schw}.

From the mid 1990's, it also became apparent that Colombeau algebras could be a significant tool with 
which to study singular problems in various geometrical settings. In particular, early work centered on 
applications to problems in General Relativity (see \cite{SVsurvey} for a recent survey), 
and Lie group analysis of partial differential equations (e.g., \cite{symm,book}).

At the same time, structural properties of Colombeau algebras came to the fore
in the work of several research groups. In particular, a thorough study of algebraic 
properties was carried out (e.g., \cite{A, V1}) and topological and functional analytic 
structures on Colombeau spaces were developed and refined to a high degree (e.g.,
\cite{S0, S, G1, G2}).

The aim of this contribution is to provide an overview of some of these developments
that show significant potential both for the intrinsic understanding of algebras of generalized
functions and for applications in geometry, differential equations, and mathematical physics.

\section*{2. Non-smooth differential geometry.} 
Throughout this paper we will employ the so-called special (or simplified) version of
Colombeau's algebras. To fix notations we briefly recall the definition of the
Colombeau algebra $\G(M)$ on a manifold $M$ (see, e.g., \cite{book}). Let ${\mathcal P}(M)$
denote the space of linear differential operators on $M$. $\G(M)$ is defined as the quotient
space $\esm(M)/\ns(M)$, where the spaces of moderate resp.\ negligible nets are defined by
\beas
\esm(M) &=& \{(u_\eps)_\eps\in \cinfty(M)^{(0,1]} \ :\ \forall K \comp M\,
\,\forall P\in {\mathcal P}(M)\, 
\, \exists\, l \\
&& \ \sup_{x\in K} |Pu_\eps(x)| = O(\eps^{-l})\}\\
\ns(M) &=& 
\{(u_\eps)_\eps\in \esm(M)\qquad :  \ \forall K\comp M \,\,\, 
\, \ \ \ \, \forall m\, 
\sup_{x\in K} |u_\eps(x)| \ \,\, = O(\eps^{m})\}
\eeas
Here and in what follows we will assume that all representatives of generalized
functions in fact depend smoothly on the regularization parameter $\eps$.
A similar definition can be given for the space $\Gamma_\G(M,E)$ of generalized sections
of a vector bundle $E\to M$, and we have the fundamental $\cinfty(M)$-module isomorphism
$$
\Gamma_{\G}(M,E) \cong \G(M) \otimes_{\cinfty(M)} \Gamma(M,E),
$$
i.e., generalized sections may be viewed globally as sections with generalized coefficient
functions. Based on regularization operations via convolution in charts (cf.\ the de Rham
regularizations in \cite{deR}) it can be shown that there exist injective sheaf morphisms
$$\iota: \Gamma(\_\,,E)
\hookrightarrow {\mathcal D}'(\_\,,E) \hookrightarrow \Gamma_{\G}(\_\,,E).
$$
An important feature distinguishing Colombeau generalized functions from Schwartz 
distributions is the availability of a point value characterization: we call a net
$(x_\eps)_\eps$ of points in $M$ compactly supported if $x_\eps$ remains in some
compact set for $\eps$ small. Two compactly supported nets $(x_\eps)_\eps$, 
$(y_\eps)_\eps$ are called equivalent, $(x_\eps)_\eps \sim (y_\eps)_\eps$, if 
$d_h(x_\eps,y_\eps)=O(\eps^m)$ $\forall m$, where $d_h$ is the distance function
induced by any Riemannian metric $h$ on $M$. The quotient space $\tilde M_c:= M^{(0,1]}$ 
is called the space of compactly supported generalized points. Then we have
(cf.\ \cite{book}, Th.\ 3.2.8):

\noindent{\bf Theorem 1.}{\em Let $u\in \G(M)$. Then $u=0$ if and only if $u(\tilde x) = 0$
for all $\tilde x \in \tilde M_c$.}

As a first pointer at algebraic properties of $\G$, let us have a look at 
the question of (multiplicative) invertibility in both $\G(M)$ and
the ring of constants in $\G(M)$ (or space of generalized numbers),
$\tilde \K$.\medskip\\

\noindent{\bf Lemma 1. }{\em Let $u\in \G(M)$. The following are equivalent:
\begin{itemize}
\item[(i)] $u$ is invertible.
\item[(ii)] $u(\tilde x)$ is invertible in $\tilde \K$ for all $\tilde x\in \tilde M_c$.
\item[(iii)] $u$ is {\em strictly nonzero}, i.e., $\forall K\comp M$ $\exists q$ s.t.\
$\inf_{p\in K}|u_\eps(p)| > \eps^q$ for $\eps$ small.
\end{itemize}
}
Similarly, for generalized numbers we have:\medskip\\

\noindent{\bf Lemma 2. }{\em Let $r\in \tilde \K$. The following are equivalent:
\begin{itemize}
\item[(i)] $r$ is invertible.
\item[(ii)] $r$ is not a zero divisor.
\item[(iii)] $r$ is strictly nonzero.
\item[(iv)] For every representative $(r_\eps)_\eps$ of $r$ there exists some 
$\eps_0>0$ such that $r_\eps \not=0$ for all $\eps<\eps_0$.
\end{itemize}
}
While the other conditions in Lemmas 1 and 2 are well-known (cf.\ \cite{book}),
(iv) from Lemma 2 is a rather recent and very convenient observation from \cite{M}.

For applications in general relativity, a notion of generalized (pseudo-)Rie\-mann\-ian 
metric is of central importance. Denoting $\Gamma_\G(M,T^0_2M)$ by $\G^0_2(M)$
we have the following characterization (\cite{gprg, M}): \medskip\\

\noindent{\bf Theorem 2. }{\em
Let $g\in \G^0_2(M)$. The following are equivalent:
\begin{itemize}
\item[(i)] $g: {\G}^1_0(M)\times {\G}^1_0(M) \rightarrow \G(M)$ is symmetric and
$\det(g)$  is invertible in $\G(M)$.
\item[(ii)] For each chart $(\psi,V)$, $\forall \tilde x \in (\psi(V))^\sim_c$: 
$\psi_*g(\tilde x)$: $\tilde \K^n \times \tilde \K^n \to \tilde \K$ is symmetric and nondegenerate.
\item[(iii)] $\det(g)$ is invertible in $\G(M)$ and $\forall \overline{V}\comp M$ 
there exists a representative $(g_\eps)_\eps$, such that each $g_\eps|_{V}$ 
is a smooth pseudo-Riemannian metric. 
\end{itemize} 
Moreover, if $g$ satisfies these equivalent conditions then $g$ has index $j$
if and only if for each chart $\psi$ and each $\tilde x$, $\psi_*g(\tilde x)$ is a symmetric
bilinear form on $\tilde \R^n$ with index $j$.
}
As in the smooth setting, the following fundamental lemma shows that each
generalized pseudo-Riemannian metric induces a unique Levi-Civita connection (\cite{gprg}):\medskip\\

\noindent{\bf Theorem 3. }{\em For any generalized pseudo-Riemannian metric $g$ on 
$M$ there exists a unique connection $\hat \nabla: {\gs}^1_0(M)\times{\gs}^1_0(M)\to{\gs}^1_0(M)$ 
such that:
\begin{itemize}
\item[\hspace*{1cm}($\nabla1$)] $\hat \nabla_X Y$  is ${\tilde \R}$-linear in $Y$.
\item[($\nabla2$)] $\hat \nabla_X Y$ is $\gs(M)$-linear in $X$.
\item[($\nabla3$)] $\hat \nabla_X(uY)=u\,\hat \nabla_X Y+X(u)Y$ for all $u\in\gs(M)$.
\item[($\nabla4$)] $[X,Y]=\hat \nabla_X Y-\hat \nabla_YX$ 
\item[($\nabla5$)] $X\langle Y,Z\rangle=\langle 
\hat \nabla_X Y,Z\rangle+\langle Y,\hat \nabla_X Z\rangle$
\end{itemize}
}
With these tools at hand, one can proceed to analyzing curvature quantities
and geodesics for singular metrics. We refer to \cite{SVsurvey} for a recent
overview of applications in general relativity. More generally, generalized
connections in principal fiber bundles have been studied in \cite{connections}.
Notions like curvature, holonomy and characteristic classes can then be
modelled in a non-smooth setting. First applications to singular Yang-Mills
equations can also be found in \cite{connections}.

A further aspect of Colombeau algebras that allows to go beyond the distributional
setting is the notion of generalized functions taking values in differentiable
manifolds (\cite{gfvm,gfvm2}). The basic idea is to consider, for given manifolds
$M$, $N$, a quotient construction on subspaces of $\mathcal{E}[M,N]:=\cinfty(M,N)^{(0,1]}$.
The corresponding growth conditions can either be modelled by asymptotic estimates
in charts (\cite{gfvm}) or, more elegantly, using `referee functions' for testing
for moderateness resp.\ negligibility, as follows: we say that a net $(u_\eps)_\eps \in
\mathcal{E}[M,N]$ (depending smoothly on $\eps$) is c-bounded if for each $K\comp M$ 
there exists some $K'\comp N$
and some $\eps_0$ such that $u_\eps(K)\subseteq K'$ for all $\eps <\eps_0$.
A c-bounded net $(u_\eps)_\eps$ is called moderate if 
$(f\circ u_\eps)_\eps\in\esm(M)\ \forall\ f\in\cinfty(N)$.
The space of moderate nets is denoted by $\esm[M,N]$. Two elements $(u_\eps)_\eps$,
$(v_\eps)_\eps$ of $\esm[M,N]$ are called equivalent, $(u_\eps)_\eps\sim (v_\eps)_\eps$,
if $(f\circ u_\eps-f\circ v_\eps)_\eps\in\ns(M)\ \forall\ f\in\cinfty(N)$. The space
of Colombeau generalized functions on $M$ taking values in $N$ is then given by
$\G[M,N]:=\esm[M,N]/\sim$.

Manifold-valued generalized functions are a necessary prerequisite for addressing
problems like determining geodesics of singular metrics or flows of generalized vector
fields. Based on $\G[M,N]$, a functorial theory of manifold-valued generalized
functions and generalized vector bundle homomorphism has been developed in 
\cite{gfvm,gfvm2}. 

On the structural level, a basic question is whether $\G[\_,N]$ forms a sheaf.
Due to the lack of algebraic structure on the target space $N$, the usual 
tools like partitions of unity are not directly available to answer this question.
Nevertheless, we have:\medskip\\

\noindent{\bf Theorem 4. }{\em $\G[\_,N]$ is a sheaf of sets.
}

\noindent{\em Sketch of proof.} The nontrivial part is to show that any coherent
family of locally defined generalized maps is given as a family of
restrictions of one globally defined generalized map.
The strategy is to use a Whitney embedding of $N$
into some $\R^n$ and then apply a gluing procedure in $\R^n$ based on partitions
of unity. In order to obtain a global representative taking values in $N$, the
retraction map of a tubular neighborhood of $N$ in $\R^n$ is employed. 
For details, see \cite{sheaf}. $\Box$

By a similar method, we obtain the following result on the inclusion of 
continuous maps in $\G[M,N]$ ($\sigma$ denotes the identical embedding
$f\mapsto (f)_\eps$ of $\cinfty(M,N)$ in $\G[M,N]$), cf.\ \cite{sheaf}:\medskip\\

\noindent{\bf Theorem 5. }{\em 
There exists an embedding $\iota: {\mathcal C}(M,N) \hookrightarrow \G[M,N]$ with the
following properties:
\begin{itemize}
\item[(i)] $\iota$ is a sheaf morphism.
\item[(ii)] $\iota|_{\cinfty(M,N)} = \sigma$.
\item[(iii)] $\iota(u)_\eps$ converges to $u$ uniformly on compact sets.
\end{itemize}
}
As an added benefit, the construction of $\G[M,N]$ provides a blueprint for defining
a space of manifold-valued distributions, as follows: set
${\mathcal A}[M,N] = \{u\in \gs[M,N] \mid \forall f\in \cinfty(N), \exists \lim_{\eps\to 0}f\circ u_\eps
\in {\mathcal D}'\}$ and let $u\approx_{\mathcal M} v$ if for all $f\in \cinfty(N)$, 
$f\circ u_\eps - f\circ v_\eps \to 0$ in ${\mathcal D}'$.
Then set ${\mathcal D}'(M,N):= {\mathcal A}[M,N]/\approx_{\mathcal M}$. For $M$,
$N$ Euclidean spaces, ${\mathcal D}'(M,N)$ singles out a subspace of bounded
distributions. Further properties (e.g., the relationship to Young measures)
are analyzed in \cite{sheaf}. 

\section*{3. Some algebraic aspects of Colombeau algebras on manifolds.}
Based on the construction in the previous sections, here we will give a few 
examples indicating the increasingly important role that an understanding 
of the algebraic structure of Colombeau-type spaces plays in a geometrical context.

To begin with, let us consider the structure of the
space of algebra isomorphisms from $\G(M)$ to $\G(N)$. In the smooth setting, it has been
known for a long time that for any algebra isomorphism $\phi: \cinfty(M) \to \cinfty(N)$
there is a unique diffeomorphism $f:N\to M$ of the underlying manifolds such that $\phi$
is given as the pullback map under $f$: $\phi = u \mapsto u\circ f$. The analogous 
problem for isomorphisms of Colombeau algebras has only recently been solved by 
H.\ Vernaeve in \cite{V2}. The result is based on the solution of `Milnor's exercise'
in the Colombeau setting, i.e.\ the characterization of multiplicative linear functionals on $\G$:\medskip\\

\noindent{\bf Theorem 6. }{\em
Every multiplicative linear functional on $\gs(M)$ is of the form
$$
e\delta_{\tilde x}: u \mapsto e u(\tilde x)
$$
for $\tilde x$ a generalized point and $e\in \tilde \K$ idempotent.
}
Using this result, we obtain\medskip\\

\noindent{\bf Theorem 7. }{\em
Let $\phi: \gs(M)\to \gs(N)$ be an algebra-isomorphism $($with $\phi(1)=1)$.
Then $\phi = f^*$ for some $f\in \gs[N,M]$ such that $f^{-1}\in \gs[M,N]$. Also, $\phi^{-1} = f_*$.
}
Next, let us investigate generalized de Rham cohomology. By 
$\Omega^p_\gs(M) = \Gamma_\gs(M,$ $\Lambda^p(M))$ we denote the space of generalized $p$-forms on $M$.
Also, as in the smooth setting we introduce the cohomology spaces by
\begin{eqnarray*}
Z^p_\gs(M) &:=& \{\omega\in \Omega^p_\gs(M) \mid d\omega = 0\} \\
B^p_\gs(M) &:=& \{\omega\in \Omega^p_\gs(M) \mid \exists \tau\in \Omega^{p-1}_\gs: \omega=d\tau\}\\
H^p_\gs(M) &:=& Z^p_\gs(M)/H^p_\gs(M)
\end{eqnarray*}
The relationship between generalized and smooth de Rham cohomology is as follows:\medskip\\

\noindent{\bf Theorem 8. }{\em 
For any $p\ge 0$ we have the following isomorphism of real vector spaces:
$$
H^p_\gs(M) \cong \tilde \R \otimes_\R H^p(M)
$$
}

\noindent{\em Sketch of proof.} Both
$$
0 \longrightarrow \ker(d) \stackrel{d}{\longrightarrow} \Omega_\gs^0(M)
\stackrel{d}{\longrightarrow} \Omega_\gs^1(M) \stackrel{d}{\longrightarrow}\dots
$$  
and
$$
0 \longrightarrow \ker(d) \stackrel{id\otimes d}{\longrightarrow} \tilde \R \otimes_\R C^\infty(M,\R)
\stackrel{id\otimes d}{\longrightarrow}  \tilde \R \otimes_\R \Omega^1(M)
\stackrel{id\otimes d}{\longrightarrow}\dots
$$
are fine resolutions of the sheaf of locally constant Colombeau generalized functions.
The result therefore follows from the abstract de Rham theorem. For details, see
\cite{connections}. $\Box$

This means that the structural difference between generalized and smooth de Rham 
cohomology is encoded precisely in the algebraic structure of the ring of
generalized numbers. 

Finally, let us return to the algebraic foundations of pseudo-Riemannian geometry in
the Colombeau setting. As can be seen from Th.\ 2, the study of bilinear forms on
$\tilde \R^n$ is of central importance here. We have (\cite{M}):\medskip\\

\noindent{\bf Theorem 9. }{\em
Let $v\in \tilde\R^n$. The following are equivalent:\\
(i) For any positive definite bilinear form $h$, $h(v,v)>0$.\\
(ii) $v$ is free (i.e., $\lambda v = 0 \Rightarrow v=0$).\\
(iii) $v$ can be extended to a basis of $\tilde\R^n$.\\
(iv) For each representative $(v_\eps)_\eps$ there exists some $\eps_0$ such that 
for all $\eps<\eps_0$, $v_\eps\not=0$.
}

Based on this result, causality notions (time-like, space-like, and null vectors) can be 
introduced and analyzed in the generalized setting. Applications include energy methods
for solving wave equations on singular space-times (cf.\ \cite{GMS}).
\section*{4. Algebraic properties of $\tilde \K$.} In this section we give a brief 
overview of known results on the algebraic structure of the ring of generalized 
numbers. For details and proofs we refer to the original sources \cite{A, AJOS, V1}.
In what follows, topological properties always refer to the sharp topology on
$\tilde \K$ (cf.\ the following section).
\begin{itemize}
\item $\tilde \K$ is a reduced ring, i.e., there are no nontrivial nilpotent elements.
\item Elements of $\tilde \K$ are either invertible or zero-divisors (cf.\ Lemma 2).
\item $e\in \tilde \K$ is idempotent ($e^2 = 1$) iff $e=e_S$, the characteristic
function of some $S\subseteq (0,1]$. 
\item $\tilde \K$ possesses uncountably many maximal ideals.
\item $\tilde \K$ is a complete topological ring.
\item The closure of any prime ideal is maximal. Conversely, every maximal ideal is closed.
\item Let $I$ be an ideal in $\tilde \K$. Then the closure of $I$ is the intersection
of all maximal ideals containing $I$.
\item $\tilde \K$ is {\em not:}
\begin{itemize}
\item Artinian
\item Noetherian
\item von Neumann regular
\end{itemize}
\item Every ideal $I$ in $\tilde \K$ is convex ($x\in I$, $|y|\le |x|$ $\Rightarrow$ $y\in I$).
\item An ideal $I$ is prime iff it is pseudoprime and radical, i.e.:
\begin{itemize}
\item $\forall S\subset (0,1]$: $e_S\in I$ or $e_{S^c}\in I$, and
\item $\forall x\in I$: $\sqrt{|x|}\in I$.
\end{itemize}
\end{itemize}
We note that many of the corresponding properties for $\G$ instead of $\tilde\K$
are the subject of ongoing research. We conclude this section with the following
interesting connection to the nonstandard space of
asymptotic numbers (cf.\ \cite{OT}), established in \cite{V1}, Th.\ 7.2:\medskip\\

\noindent{\bf Theorem 10. }{\em
Let $I$ be a maximal ideal in $\tilde \K$. Let $\mathcal{U}:=\{S\subseteq (0,1]\mid e_{S^c}\in I\}$.
Let ${}^*\K$ be the nonstandard field constructed by the ultrafilter $\mathcal U$ and
let $\rho$ be the infinitesimal with representative $(\eps)_\eps$. Then ${}^\rho\K$ is
canonically isomorphic to $\tilde\K/I$.
}
\section*{5. Topology and functional analysis.} Topologies on spaces of Colombeau
generalized functions and generalized numbers were originally introduced by
D.\ Scarpalezos by the name of {\em sharp topologies} in 1993 (and published only 
later in \cite{S0,S}). After the field lay dormant for some years (in which the
main focus of research was on applications in PDEs) there occurred a veritable surge 
of activities lately. In particular, the fundamental work by C.\ Garetto \cite{G1,G2}
has led to the development of a full-scale locally convex theory for algebras of 
generalized functions. In this section we outline some of the main features of this
theory.

For any given locally convex vector space $E$ whose topology is induced by the family of
seminorms $(p_i)_{i\in I}$, we set
\begin{eqnarray*}
\mathcal{M}_E &:=& \{(u_\eps)_\eps \in E^{(0,1]} \mid \forall i\, \exists N\,: p_i(u_\eps)= O(\eps^{-N})\}\\
\mathcal{N}_E &:=& \{(u_\eps)_\eps \in E^{(0,1]} \mid \forall i\, \forall q\,: p_i(u_\eps)= O(\eps^{q})\}\\
\gs_E &:=& \mathcal{M}_E/\mathcal{N}_E
\end{eqnarray*}
Then $\gs_E$ is a $\tilde \C$-module. The special Colombeau algebra
$\gs(\Omega)$ is obtained as the special case $E=\cinfty(\Omega)$
of this construction (cf.\ \cite{G1,DHPV}). 

On $\G_E$ we introduce valuations given by 
$$
v_{p_i}(u):= \sup\{b\in \R \mid p_i(u_\eps)=O(\eps^b)\}
$$
(here $(u_\eps)_\eps$ is any representative of $u\in \G_E$). The valuations, in turn,
induce ultra-pseudo-seminorms (ups) via 
$$
\mathcal{P}_i := e^{-v_{p_i}}.
$$
This family of ups defines the sharp topology on $\G_E$. 
As an important special case we may take $E = \C$, in which case $\gs_E = \tilde \C$.
Here we only have one seminorm, $p(x) = |x|$ which induces a valuation $v$ and 
a corresponding ups denoted by $|\,\,|_e$.

More generally we may introduce suitable notions for directly generalizing 
locally convex topologies to the $\tilde \C$-module setting. Recall that
for $V$ a vector space and $X\subseteq V$, $X$ is called {\em absorbent} in $V$ if 
$\forall u\in V$ $\exists \lambda_0$ $\forall \lambda\ge \lambda_0$:
$u\in \lambda X$. Let now $\gs$ be a $\tilde \C$-module and let  
$A\subseteq \gs$. If we let $\lambda_0$ correspond to the infinitesimal
$[(\eps^a)_\eps]$, then $\lambda_0 \cong [(\eps^a)_\eps] \le [(\eps^b)_\eps]\cong 
\lambda$ iff $b\le a$, so we are led to defining: 
$A$ is called $\tilde \C$-{\em absorbent} if $\forall u\in \gs$ $\exists a\in \R$ 
$\forall b\le a$: $u\in [(\eps^b)_\eps]A$. Similarly, we call 
$A$ $\tilde \C$-{\em balanced} if $\forall \lambda\in \tilde \C$ with $|\lambda|_e\le 1$:
$\lambda A \subseteq A$.

To introduce a suitable notion of convexity, recall that a subset $X$ of a vector space
$V$ is a convex cone in $V$ if $X+X\subseteq X$ and $\forall \lambda\in (0,1]$: 
$\lambda X\subseteq X$. Thus we call a subset $A$ of a $\tilde \C$-module $\G$
$\tilde \C$-{\em convex} if $A+A\subseteq A$ and $\forall b\ge 0$:
$[(\eps^b)_\eps]A\subseteq A$. Finally, we define a 
locally convex topological $\tilde \C$-module to be a toplogical $\tilde \C$-module
(which means that $+$ and $\lambda \cdot$ are continuous) with a base of 
$\tilde \C$-convex neighborhoods of $0$.

This provides the starting point for a by now highly developed theory of locally
convex $\tilde \C$-modules which to a large extent parallels the theory of locally
convex vector spaces. Some of the main features of the theory are:
\begin{itemize}
\item The ups take over the role of seminorms.
\item Completeness, metrizability, projective and inductive limits have been studied.
\item There is a theory of barrelled and bornological $\tilde \C$-modules.
\item Examples: $\gs_c(\Omega)$ (corresponding to $\mathcal{D}(\Omega)$) is a strict inductive limit.
$\gs(\Omega)$ (corresponding to $\cinfty(\Omega)$) is a Frechet $\tilde \C$-module. The standard
spaces $\gs_\tau(\Omega)$, $\gs_{\mathcal{S}}(\Omega)$, $\gs^\infty(\Omega)$, etc.\
can all be treated within the theory.
\item Duality theory, study of 
$$
\mathcal{L}(\gs,\tilde \C):= \{T:\gs\to \tilde \C \mid T\ \tilde\C\mathrm{-linear\ and\ continuous}\}
$$
An example is the generalized delta distribution (point evaluation at
$\tilde x\in \tilde \Omega$): $\delta_{\tilde x} = u\mapsto u(\tilde x) \in 
\mathcal{L}(\gs_c(\Omega),\tilde \C)$.
\item Based on this, kernels of pseudodifferential operators can be constructed as elements of 
$\mathcal{L}(\gs_c(\Omega\times\Omega),\tilde \C)$ (cf.\ also \cite{D}).
\item Microlocal analysis in the dual of Colombeau algebras, see \cite{G3}.
\item A Hahn-Banach theorem is not attainable in general due to algebraic
obstructions (\cite{V1}). 
\item Several open mapping and closed graph theorems 
and applications to $\gs^\infty$-hypo\-ellip\-ticity are given in \cite{G4}.
\end{itemize}
\section*{6. Conclusions and Outlook.} As can be seen from the above summary
of results, Colombeau theory is currently undergoing a profound and
far-reaching conceptual restructuring. Several branches of research that
so far had been rather disconnected have seen fruitful and promising 
interactions. As a first example we have seen the strong links between global
analysis and algebraic properties in Section 3. These will give rise to 
a new {\em algebraic} approach to non-smooth differential geometry.
Moreover, the algebraic causality structures also mentioned in Section 3
are currently being pursued as a tool for generalizing the Hawking and Penrose
singularity theorems of general relativity to space-times of low differentiability. 

Interactions between algebra and PDE theory include topics like a refined
study of hypoellipticity properties (which will require at least the
rudiments of real algebraic geometry in the generalized setting). Moreover,
as was indicated in Section 5, there are by now strong ties between functional
analytic methods and the theory of pseudodifferential and Fourier integral 
operators. Similarly, there are close connections between such methods
and variational problems of low regularity.

There already are examples of abstract (functional analytic)
existence results for concrete analytical problems in PDE theory in the Colombeau
framework, a direction of research which without doubt will gain
importance in the near future. The hope here is to provide a toolkit
(similar to the one available in classical analysis) of topological
and algebraical methods for solving problems of non-smooth analysis
and geometry.


\end{document}